\documentclass[11pt,twoside]{article}
\usepackage{amssymb,amsfonts,amsmath,latexsym,graphicx}

\newtheorem{lemma}{Lemma}
\newtheorem{proposition}[lemma]{Proposition}
\newtheorem{corollary}[lemma]{Corollary}
\newtheorem{theorem}{Theorem}
\newtheorem{remark}{Remark}

\renewcommand{\appendix}[2]{\section*{#1}\label{#2}\addtocontents{toc}{\vskip 0pt {\hspace*{-12pt}\bf #1\hfill\ \thepage}}}

\newcommand{\email}[1]{{\small E-mail: {\textsf {#1}}}}
\newcommand{\http}[1]{{\small Internet: {\textsf {#1}}}}
\newcommand{\affil}[1]{{\small\sl #1}}
\newcommand{\keywords}{\medskip\noindent{\bf Keywords. }}
\newcommand{\AMS}{\noindent{\small\bf AMS classification (2000). }}
\newcommand\runninghead[2]{\pagestyle{myheadings}\markboth 
 {{\footnotesize\it{\quad #1}\hfill}}
 {{\footnotesize\it{#2\hfill\quad}}}}\headsep=40pt
\newcommand{\proof}{\noindent {\sl Proof.\/} }

\newcommand{\eqn}[1]{(\ref{#1})}
\newcommand{\be}{\begin{equation}}
\newcommand{\ee}{\end{equation}}
\newcommand{\N}{\mathbb{N}}
\newcommand{\R}{\mathbb{R}}
\newcommand{\Z}{\mathbb{Z}}
\newcommand{\ds}{\displaystyle}
\def\qed{\hfill$\diamondsuit$}
\def\D{\mathcal{D}}

\begin{document}

\runninghead{\hfill\ R. Benguria, I. Catto, J. Dolbeault, \& R. Monneau \hfill\ \today}{\today\hfill\ Oscillating minimizers \hfill\ }
\title{\sl Oscillating minimizers of a fourth order problem invariant under scaling
\footnote{\copyright\, 2003 by the authors. This paper may be reproduced, in its entirety, for non-commercial purposes.}\kern 10pt\footnote{Supported by ECOS-CONYCIT under contract \# C02E06}}
\author{
Rafael D. Benguria\footnote{\, Partially supported by Fondecyt (Chile), 
project \# 102-0844.}\\
\affil{Departamento de F\'\i sica P. U. Cat\'olica de Chile, Casilla 
306, Santiago 22, Chile}\\
\affil{Phone: (56) 2 354 4486, Fax: (56) 2 553 6468}\\
\email{rbenguri@fis.puc.cl},
\http{http://www.fis.puc.cl/\raisebox{-4pt}{$\widetilde{\phantom{x}}$}rbenguri/}\\
\\
Isabelle Catto, Jean Dolbeault\footnote{\, Partially supported by the european networks ``Analysis and Quantum'' and ``HYKE'' under contracts \# HPRN-CT-2002-00277 and \# HPRN-CT-2002-00282}
 \\
\affil{Ceremade (UMR CNRS no. 7534), Universit\'e Paris IX-Dauphine,}\\
\affil{Place de Lattre de Tassigny, 75775 Paris C\'edex~16, France}\\
\affil{Phone: (33) 1 44 05 44 65, 46 78, Fax: (33) 1 44 05 45 99}\\
\email{catto, dolbeaul@ceremade.dauphine.fr}\\
\http{http://www.ceremade.dauphine.fr/\raisebox{-4pt}{$\widetilde{\phantom{x}}$}catto, dolbeaul/}\\
\\
R\'egis Monneau\\
\affil{CERMICS, ENPC, 6 \& 8 avenue Blaise Pascal,}\\
\affil{Cit\'e Descartes, Champs-sur-Marne, 77455 Marne-La-Vall\'ee, France}\\
\affil{Phone: 01 64 15 35 80 \& 01 45 85 19 44 - Fax: 01 64 15 35 86}\\
\email{monneau@cermics.enpc.fr}\\
\http{http://cermics.enpc.fr/\raisebox{-4pt}{$\widetilde{\phantom{x}}$}monneau/home.html}\\
\date{}
}
\maketitle
\thispagestyle{empty}\vskip -0.9cm
\noindent {\bf Abstract.} {\sl By variational methods, we prove the inequality\\
\centerline{$\ds\int_\R u''{}^2\,dx-\int_\R u''\,u^2\,dx\geq I\,\int_\R u^4\,dx\quad \forall\; u\in L^4(\R)\;\mbox{such that}\; u''\in L^2(\R)$}
for some constant $I\in (-9/64,-1/4)$. This inequality is connected to Lieb-Thirring type problems and has interesting scaling properties. The best constant is achieved by sign changing minimizers of a problem on periodic functions, but
does not depend on the period. Moreover, we completely characterize the minimizers of the periodic
problem.}\medskip

\medskip

\AMS{Primary: 35J35, 26D20 -- Secondary: 47J20, 49J40}
%
%
%
%
%
%

%
%
%

\section{Introduction}\label{Intro}

In this paper, we prove the inequality 
\begin{equation}
\label{MainIneq} 
\int_\R u''{}^2\,dx-\int_\R u''\,u^2\,dx\geq I\,\int_\R u^4\,dx\quad \forall\; u\in
L^4(\R)\;\mbox{such that}\; u''\in L^2(\R)
\end{equation}
for some negative constant $I$. This inequality is a special case of more general inequalities involving terms like~: $u''{}^2$, $u''\,u^2$, $u'{}^2\,u$, $u'{}^4/u^2$, $u^4$... which all share the same scaling behaviour under the scaling
$\sigma\mapsto \sigma^2\,u(\sigma\,\cdot)$.
\medskip

Apart from its own interest, the initial motivation for studying such a problem is connected with Lieb-Thirring inequalities. In \cite{Benguria-Loss00}, R.D. Benguria and M. Loss gave a simple proof of a theorem of A. Laptev and T. Weidl \cite{Laptev-Weidl00} using a commutator method. It was then natural to ask if such a method could also work for fourth order operators as well. This has been recently investigated by A. Laptev and J. Hoppe \cite{Laptev-Hoppe03}. It turns out that the above inequality plays an important role for such an approach.\medskip

Our main result is the following.
\begin{theorem}\label{Main} The best constant $I$ in Inequality \eqn{MainIneq} is given by
\[I=\inf\left\{\frac{\int_0^T u''{}^2\, dx-\int_0^T u''\,u^2\,dx}{\int_0^T u^4\,dx}\;:\; u\not\equiv 0,\,u\in C^\infty(\R/T\Z)\right\}\] 
where $C^\infty(\R/T\Z)$ denotes the set of $T$-periodic functions in $C^\infty(\R)$. The best constant is not achieved on $\R$ but it is achieved on the set of periodic functions, and it is independent of the period $T$. It takes values in $(-1/4,-9/64)$.

Moreover, for any $T>0$, there exists a unique minimizer with minimal period $T$, up to translations. This minimizer changes sign.\end{theorem}
The difficulty of the above minimization problem comes from the loss of compactness due to 
the scaling and translation invariances. It is furthermore interesting to understand the rather
nonstandard properties of the minimizers in the periodic case, which for instance are always
given by sign changing functions. On the whole real line, we will show that minimizing sequences
can be chosen as the restriction to a finite number of periods of periodic functions, up to some
tail, whatever the period is, and that the infimum is reached when the number of periods goes to
infinity.\medskip

A result similar to Theorem \ref{Main} was obtained by A. Leizarowitz and V.J. Mizel~\cite{LM} for some infinite-horizon variational problems of second order leading to a fourth order ODE. Certain conditions where given in \cite{MZ} to assure the uniqueness (up to translation) of the periodic minimizer. For a similar ODE, V.J. Mizel, L.A. Peletier and W.C. Troy proved in \cite{MPT} (also see \cite{VanDenBerg-Peletier01}) that any periodic minimizer has to be even with respect to its extrema and is therefore a single-bump function. L.A. Peletier in \cite{P} proved using a cut-and-paste argument in the $(u,u')$-plane that the map $x\mapsto (u(x),u'(x))\in \R^2$ is injective. However, the specificity of the problem considered in this paper is the scaling invariance which is not present in the above mentioned references.\medskip 

This paper is organized as follows. We first state some preliminary results in
Section~\ref{Prelim}. Then we prove Theorem~\ref{Main} and some qualitative properties of the minimizers in Section~\ref{Proof}. The last Section is devoted to numerical computations of the best constant, whose
value~is
\[ I=-0.1580...\]
and for which precise theoretical estimate still need to be found.

\section{Preliminary results}\label{Prelim}

Let us define
\begin{equation}\label{eq:def-I}
I:=\inf\Big\{Q_\R(u)\;:\; u\not\equiv 0\,,\; u\in L^4(\R)\,,\; u''\in L^2(\R)\Big\}
\end{equation}
where
\begin{equation}\label{eq:def-Q}
Q_\R(u):=\frac{\int_\R u''{}^2\,dx-\int_\R u''\,u^2\,dx}{\int_\R u^4\,dx}\;.
\end{equation}
By a density argument,
\begin{equation}\label{eq:def-I2}
I:=\inf\{Q_\R(u)\;:\;u\not\equiv 0\,,\; u\in \D (\R)\}\;.
\end{equation}
The analogous variational problem with periodic boundary conditions on $[-T,T)$ reads
\begin{equation}\label{eq:Iper}
I_T:=\inf\Big\{\,Q_T(u)\;:\;u\not\equiv 0\,,\; u\in L_{\rm loc}^4(\R)\,,\; u''\in L_{\rm loc}^2(\R)\,,\; u(\cdot+2\,T)=u\,\Big\}\;,
\end{equation}
where
\begin{equation}\label{eq:Qper}
Q_T(u):=\frac{\int_{-T}^{+T} u''{}^2\,dx-\int_{-T}^{+T} u''\,u^2\,dx}{\int_{-T}^{+T} u^4\,dx}\;.
\end{equation}
In the rest of this paper, we prefer to work with $2\,T$ periodic functions instead of~$T$ periodic functions and consider $[0,T)$ 
as the standard half period, for notational convenience. We also denote more generally
\[
Q_J(u):=\frac{\int_J u''{}^2\,dx-\int_J u''\,u^2\,dx}{\int_J u^4\,dx}\;,
\]
for any interval $J$ of $\R$. $Q_\R$ will sometimes be simply denoted by $Q$ 
when there is no ambiguity. We shall prove in the following that $I_T=I$ for any $T>0$ and then prove a series of results on the features of the minimizers.
\begin{lemma}{\rm [Well definiteness - First rough estimates]}
$$-\frac{1}{4}\leq I < 0\;.$$
\end{lemma}
\proof Let $u\in\D (\R)$ with $u\leq 0$, $u\not \equiv 0$. We observe that for smooth enough functions
\[ 
-\int_\R u''\,u^2\,dx=2\,\int_\R u\,u'{}^2\,dx
\]
by integrating by parts. Then, for every $\lambda>0$,
\[
Q(\lambda\, u)=\lambda^{-2}\,\frac{\int_\R u''{}^2\,dx}{\int_\R u^4\,dx}+2\,\lambda^{-1}\,\frac{\int_\R u\,u'{}^2\,dx}{\int_\R u^4\,dx}\,,
\]
and the second term, which can be taken negative by choosing $u$ non-positive, dominates as $\lambda$ goes to infinity. This proves the negative upper bound.

To get the lower bound, we simply observe that
\be\label{Cas14}
\int_\R u''{}^2\,dx-\int_\R u''\,u^2\,dx=\int_\R \Big\vert
u''-\frac 12\,{u^2} \Big\vert^2\,dx-\frac 14\,\int_\R u^4\,dx\;.
\ee
\qed
\vskip6pt
As claimed in the introduction, the variational problem has some scaling invariance (apart from the obvious translations invariance), which play an important role in the analysis of the minimizing sequences and their possible loss of compactness. 
\begin{lemma}{\rm [Scaling invariance]}\label{lem:scaling}
For every $u\not\equiv 0$ such that $u\in L^4(\R)$, $u''\in L^2(\R)$ and for every $\sigma>0$, if we define
$u_\sigma:=\sigma^2\,u(\sigma\cdot )$, then
\begin{equation}
\label{eq:scaling} Q(u_\sigma)=Q(u)\;.
\end{equation}
Similarly, for any $u\in L^4(0,T)$ such that $u''\in L^2(0,T)$, 
\[
Q_T(u_\sigma)=Q_{\sigma T}(u)\;.
\]
Therefore, for every $T>0$, $I_T=I_1$.
\end{lemma}
\vskip6pt
The proof is straightforward and left to the reader. We now prove that the variational problem over $\R$ reduces to the same problem but stated on periodic functions. 
\begin{lemma}{\rm [Reduction to periodic functions]}\label{lem:reduc} For any $T>0$,
\[I=I_T\;.\]
\end{lemma}
\proof Let $\varepsilon>0$ and let $u\in \D (\R)$, $u\not \equiv 0$, be
such that $I\leq Q(u)\leq I+\varepsilon$. For $T$ large enough so that ${\rm supp}(u)\subset [-T,T]$, $u$ may be replicated as a
$C^\infty$ periodic function and therefore $I_1=I_T\leq Q(u)\leq
I+\varepsilon$. Whence $I_1\leq I$ since $\varepsilon$ is
arbitrary.

\medskip For the reverse inequality, we argue as follows. Let $N$ be a positive integer aimed at going to infinity. Let $u_1$ be a $1$-periodic smooth function such that 
\[
Q_1(u_1)<I_1+\varepsilon\;.
\]
We may build a function $u\in H_{\rm loc} ^2(\R)$ with compact support in $[-(N+1),N+1]$ in the following way
\[
u(x)=\left\{ \begin{array}{ll}\quad 0 &\textrm{ if } |x|\geq N+1,\\ \\
\quad u_1 &\textrm{ in }[-N,N],\end{array} \right.
\]
and $u$ glues $u_1$ to 0 on $[-(N+1),-N]\cup [N,N+1]$. Then
\[
I\leq Q(u)=Q_1(u_1)+O\Big(\frac{1}{N}\Big)<I_1+O\Big(\frac{1}{N}\Big)+\varepsilon\;,\]
so that $I-I_1$ can be made arbitrarily small for $N$ large enough and $\varepsilon$ small enough.
\qed
\vskip6pt
The rest of the section is devoted to the analysis of the properties of the solutions of the associated Euler-Lagrange equations. 
\begin{lemma}{\rm [Euler-Lagrange equations and regularity]}\label{Elr} Let us assume that some function $u$ is a minimizer either of $Q$ or of $Q_T$, for some $T>0$. Then $u$ is a classical solution to the Euler-Lagrange equation
\begin{equation}\label{eq:EL}
u^{(iv)}-2\,u''\,u-u'{}^{\,2}=2\,I\,|u|^2\,u
\end{equation}
on $\R$ and $u$ is a $C^\infty$ function.\end{lemma}
\vskip6pt
Note that the Lagrange minimizer coincides with the value of the functional, which is unusual in non-linear settings. 
\vskip6pt\noindent
\proof The Euler-Lagrange equations are easily obtained by considering a variation of $Q_J$, where, here and below, $J$ stands for $\R$ or $(-T,T)$. As for the regularity, we first get for any $x$, $y\in J$
\[
\Big| u'(x)-u'(y)\Big|\leq \left|\int_x^y u''(s)\;ds\right|\leq \sqrt{|x-y|}\,\left|\int_x^y {u''}^2\;ds\right|^{1/2}
\]
by integrating between $x$ and $y$ and using the Cauchy-Schwarz inequality, so that $u$ is bounded in $C^{1,1/2}(J)$. Because of the Euler-Lagrange equation, $u^{(iv)}$ is bounded in $C^{3,1/2}(J)$ for the same reason as above. The $C^\infty$-regularity follows by bootstrapping.\qed
\vskip6pt 

This lemma now helps to better estimate the value of the infimum $I$. 

\begin{lemma}{\rm [Improved estimate]}\label{lem:improved}
$$ I < -\frac 9{64}\;.$$
\end{lemma}
\proof Let $u$ be a $C^2$ non-positive function with compact support. After one integration
by parts, we can write
\[
\int_\R u''{}^2\,dx-\int_\R u''\,u^2\,dx=-\frac 9{64}\, \int_\R u^4\,dx+\int_\R \,\Big|u''-\frac 38\,u^2-\frac 23\,\frac{{u'}^2}{u}\Big|^2\,dx\;.
\]
Let us prove first that one can find a solution to
\begin{equation}
\label{Carre}
u''-\frac 38\,u^2-\frac 23\,\frac{{u'}^2}{u}=0\;.
\end{equation}
On the support of $u$, define $y:=-|u|^{1/3}$ and solve
\[\left\{\begin{array}{l}
y''=\frac 18\,|y|^4\;,\cr\cr
y'(0)=0\;,\quad y(0)=y_0<0\;.
\end{array}\right.\]
Then 
\[
\bar u:=\left\{\begin{array}{ll}
|y|^2\,y\quad &\mbox{if}\; y<0\cr\cr
0\quad &\mbox{otherwise}
\end{array}\right.
\]
is a solution of \eqn{Carre} on the support of $\bar u$. Moreover, it is of class $C^2$ on $\R$ and 
\[
Q_\R(\bar u)=-\frac 9{64}\;.
\]
Note that on the boundary of its support, 
$\bar u'''\neq 0$. Let us extend $\bar u$ periodically. If one had $I=-9/64$, then $\bar u$ would solve the Euler-Lagrange equation \eqn{eq:EL} on $\R$ and $\bar u'''$ would be continuous, which is clearly not the case. This ends the proof.\qed
\begin{lemma}{\rm [Lower bound for $I$]}\label{lem:lower1} Let $T>0$ and assume that $Q_T$ has a nontrivial periodic minimizer with period $2\,T$. Then 
\[I_T>-\frac 1 4\;.\]
\end{lemma}
\proof If we had $I_T=-\frac 1 4$, then any minimizer would be nonpositive since
\[
\int_{u>0} u''{}^2\,dx-\int_{u>0} u''\,u^2\,dx=\int_{u>0} u''{}^2\,dx+2\,\int_{u>0} u\,u'{}^2\,dx \geq 0
\]
and because of \eqn{Cas14}, should satisfy
\[
\int_{-T}^T \Big\vert u''-\frac 12\,{u^2} \Big\vert^2\,dx =0\;.
\]
However any solution of 
\[\left
\{
\begin{array}{l}
u''=\frac 12\,{u^2}\cr\cr
u'(0)=0\;,\quad u(0)=u_0<0
\end{array}\right.\]
has a non-zero derivative at ending points $-T$ and $T$. This is again a contradiction with the regularity of
any solution of the Euler-Lagrange equation~\eqn{eq:EL}. 
\qed

\vskip6pt 
\begin{proposition}{\rm [Reduction to periodic functions that decrease on the half
period]}\label{prop:reduc} The infimum $I$ is approximated by a minimizing sequence
$(u_n)_{n\in\N}$ with the following properties: each $u_n$ has compact support and is made of the restriction to a finite number of periods of periodic sign-changing functions which are even and monotone on half of the period, up to some tail. \end{proposition}
\proof As seen in the proof of Lemma~\ref{lem:reduc}, we can choose a minimizing sequence
$(u_n)_{n\in\N}$ of $I$ as the restriction to a finite number of periods of 
periodic functions, up to some tail. The infimum is then reached when the number of periods goes
to infinity. Moreover we know from the proof of Lemma~\ref{lem:improved} that $u_n$ must be
sign-changing (at least for $n$ large enough). We thus denote by
$x_n^i$ the critical points of $u_n$ for every $n$ with $1\leq i\leq N_n$. Assume that for each
$n\in\N$, $N_n<+\infty$ and that these points are ordered~: $x_n^{i}< x_n^{i+1}$ for any~$i$. 
If $Q_{(-\infty,x_n^{i_0})}(u_n)\leq Q_{(x_n^{i_0},+\infty)}(u_n)$, we do not increase the ``energy"
$Q(u_n)$ by replacing 
$u_n(x)$ by $u_n(2x_n^{i_0}-x)$ for any $x\geq x_n^{i_0}$. Moreover, if $i_0$ is such that
$Q_{(x_n^{i_0},x_n^{i_0+1})}(u_n)\leq Q_{(x_n^{i},x_n^{i+1})}(u_n)$ for any $i\neq i_0$, it is easy to
build a new function $\tilde u_n$ which is even, $2\,\tilde T_n$-periodic on an interval
$(-N_n\,\tilde T_n,N_n\,\tilde T_n )$ with $N_n$ large and $\tilde T_n=x_n^{i_0+1}-x_n^{i_0}$, and 
such that
$\tilde u_n(x)=u_n(x+x_n^{i_0})$ for any $x\in(0,\tilde T_n)$ and
\[
Q_\R(\tilde u_n)\leq Q_\R( u_n)+O\Big(\frac{1}{N_n}\Big)
\]
(the idea is to take sufficiently many periods, {\sl i.e.\/} $N_n$ large enough, and to
then glue the function to
$0$ as in the proof of Lemma~\ref{lem:reduc}). By construction, $\tilde u_n$ is monotone on
$\big(0,\tilde T_n\big)$ and up to a shift of a half period, we may assume that it is
strictly decreasing.\qed

\section{Proof of the main result}\label{Proof}
According to the results of Section~\ref{Prelim}, the minimization problem in the whole space is
reduced to the minimization problem in the periodic case. It remains to prove that $I_T$ is achieved for some $T>0$, which is the core of the proof of Theorem \ref{Main}.
\begin{proposition}{\rm [Existence of a minimizer for the periodic case]}\label{prop:ExistPeriodic} For any
{$T>0$,} there exists a smooth nontrivial function $u$ of period $2\,T$ such that 
\[Q_T(u)=I\;.
\] 
Moreover, there is at least one minimizer $u$ which attains its absolute
maximum at $0$ (up to a translation), satisfies $u(0)>0$, $u'(0)=u'(T)=0$, and 
is even, decreasing on $(0,T)$. In addition, $u$ changes sign in $(0,T)$ and solves on $\R$ the fourth order ODE
\begin{equation}\label{eq:EL-per}
\left \{ \begin{array}{l}
u^{(iv)}-2\,u''\,u-u'{}^{\,2}=2\,I\,|u|^2\,u\;,\\ \\
u(\cdot+2\,T)=u(\cdot)\;.
\end{array}\right.
\end{equation}
\end{proposition}
\proof Let us start with some preliminary considerations. By density, the infimum of $Q$ on~$\R$ can be computed on the set of smooth functions with compact support~:
\[
I=\inf_{u\in{\cal D}(\R)}Q(u)\;.
\] 
According to Lemma \ref{lem:reduc}, it is clear that 
\[
I=\inf_{u\in C^\infty_{\scriptstyle\rm per}(\R)}Q_{T(u)}(u)\;,
\] 
where $2\,T(u)$ is the minimal period of $u$. According to Proposition~\ref{prop:reduc}, we can further
ask that $u$ is monotone decreasing on $(0,T(u))$ and even. Thus we can reduce the problem to the
case of Neumann boundary conditions
\[
I=\inf_{u\in{\cal N}}Q_{(0,T)}(u)
\] 
where $\cal N$ is the set of the $2\,T$-periodic even functions $u\in C^\infty(\R)$ such that $u'<0$ on $(0,T)$ and $u'(0)=u'(T)=0$. Because of Lemma~\ref{lem:scaling}, since for $u_\sigma:=\sigma^2\,u(\sigma\,\cdot )$,
\[
Q_{(0,T(u_\sigma))}(u_\sigma)=Q_{(0,T(u))}(u)\;,
\]
with $T(u_\sigma)=T(u)/\sigma$, there is no restriction to assume that $T(u)=1$. Thus 
\[
I=\inf_{\begin{matrix}u\in C^\infty(\R)\\
u'<0\;\mbox{\scriptsize\rm on}\;(0,1)\\
u'(0)=u'(1)=0\end{matrix}}Q_{(0,1)}(u)\;.
\] 

\medskip
To find a minimizer to the above minimization problem, we shall consider a minimizing sequence 
$(u_n)_{n\in\N}$ of the following equivalent minimization problem of Nehari type~:
\be\label{Nehari}
I=\inf_{\begin{matrix}u\in C^\infty_{\scriptstyle\rm per}(\R)\\
u'<0\;\mbox{\scriptsize\rm on}\;(0,T(u))\\
u'(0)=u'(T(u))=0\\
\int_0^{T(u)}|u|^4\, dx=1\end{matrix}}
\int_0^{T(u)}\left[u''{}^2-u''\,u^2\right]\; dx\;.
\ee
Since for $u_\mu(x):=\mu^{1/4}u(\mu\,x)$, the quantity 
\[
\mu\mapsto\int_0^{T(u_\mu)}\left[u_\mu''{}^2- u_\mu''\,u_\mu^2\right]\; dx=\mu^{7/2}\int_0^{T(u)}u''{}^2\, dx-\mu^{7/4}\int_0^{T(u)}u''\,u^2\, dx
\]
has a minimum for
\[\mu^{7/4}=\frac{\int_0^{T(u)}u''\,u^2\, dx}{2\int_0^{T(u)}u''{}^2\, dx}\;,
\]
there is no restriction to assume that the problem has already been optimized with respect to $\mu$, so that we may further impose
\[
\int_0^{T(u)}u''\,u^2\, dx=2\int_0^{T(u)}u''{}^2\, dx\;.
\]
It is then clear that 
\[ 
Q_{(0,T(u))}(u)<0
\]
since $\int_0^{T(u)}[u''{}^2-u''\,u^2]=-\int_0^{T(u)}u''{}^2$.
Going back to the minimizing sequence, we impose that
\[
\int_0^{T(u_n)}|u_n|^4\, dx=1\quad\mbox{and}\quad
\int_0^{T(u_n)}u_n''\,u_n^2\, dx=2\int_0^{T(u_n)}u_n''{}^2\, dx\;,
\]
so that 
\begin{equation}\label{eqn:limit}
\lim_{n\to\infty}\int_0^{T(u_n)}u_n''{}^2\, dx=|I|\quad\mbox{and}\quad
\lim_{n\to\infty}\int_0^{T(u_n)}u_n''\,u_n^2\, dx=2\,|I|\;.
\end{equation}

\medskip
Define now $T_n=T(u_n)$ and assume first that $\liminf_{n\to +\infty}T_n=+\infty$. 
There is no restriction to assume that $u_n$ is strictly sign changing (if not, we would obtain~:
$I\geq -9/64$)~: For any $n\in\N$, there exists an $x_n\in (0,T_n)$ such that $u_n(x_n)=0$. Let us
prove that $\liminf_{n\to +\infty}u_n'(x_n)=0$. For that purpose, consider $\bar
u_n(\cdot)=u_n(\cdot+x_n)$. Since $u_n$ is nonincreasing, 
\[u_n(x)\geq u_n(x_n-1)=\bar u_n(-1)\quad\forall\;x\in (0,x_n-1)\]
if $x_n>1$, and
\[\bar u_n(1)=u_n(x_n+1)\geq u_n(x)\quad\forall\;x\in (x_n+1,T_n)\]
if $x_n<T_n-1$. For $n$ large enough, at least one of these two conditions has to be satisfied and either 
$x_n\to +\infty$ or $T_n-x_n\to +\infty$. Since $(u_n)_{n\in\N}$ is bounded in $L^4(0,T_n)$, this means that either $\bar u_n(-1)\to 0$ or $\bar u_n(1)\to 0$. On the other hand, $(\bar u_n)_{n\in\N}$ is bounded in $H^2(-1,1)$ 
and therefore converges up to the extraction of a subsequence to a limit $\bar u$ weakly in
$H^2(-1,1)$ and strongly in $C^{1,1/2}(-1,1)$. Since $\bar u\equiv 0$ either on $(-1,0)$ or on
$(0,1)$, $\bar u'(0)=0$, this proves that $\liminf_{n\to +\infty}u_n'(x_n)=0$. 

It is then easy to check that for the minimizing sequence $(u_n)_{n\in\N}$,
we can impose $u_n'(x_n)=0$ for any $n\in\N$, up to a small change of the
sequence $(u_n)_{n\in\N}$. But this is contradictory with the fact that 
\begin{eqnarray*}
\int_0^{T_n}\left[u_n''{}^2 - u_n''\,u_n^2\right]\, dx 
&=& \int_0^{x_n}\left[u_n''{}^2 +2\,u_n\,{u_n'}^2\right]\, dx \\
&& -\frac 9{64}\int_{x_n}^{T_n}\kern -5pt |u_n|^4\,dx+\int_{x_n}^{T_n}\left|u_n''-\frac
38\,u_n^2+\frac 23\,\frac{{u_n'}^2}{u_n}\right|^2\,\kern -5pt dx \\ &&
\qquad\geq-\frac 9{64}\int_{x_n}^{T_n}|u_n|^4\,dx\;,
\end{eqnarray*}
\[-\frac 9{64}\leq Q_{[0,T_n]}(u_n^-)\leq Q_{[x_n,T_n]}(u_n)\leq Q_{[0,T_n]}(u_n)\;,\]
which proves that $(u_n^-)_{n\in\N}$ is also a minimizing sequence for $Q$ and shows that $I=-9/64$, 
a contradiction with Lemma~\ref{lem:improved}. 

\medskip Thus we know that $\limsup_{n\to +\infty}T_n<+\infty$, eventually up to the extraction of a subsequence. Let us rescale the minimizing sequence $(u_n)_{n\in\N}$~:
\begin{equation}\label{rel}
v_n(x)=T_n^{1/4}\,u_n(T_n\,x)\quad\forall\; x\in (0,1)\;,
\end{equation}
so that $v_n$ is monotone decreasing on $(0,1)$, $v_n'(0)=v_n'(1)=0$, and \eqn{eqn:limit} can be
rephrased into
\begin{equation}\label{rel2}
\lim_{n\to\infty}T_n^{-7/2}\,\int_0^1v_n''{}^2\, dx=|I|\quad\mbox{and}\quad
\lim_{n\to\infty}T_n^{-7/4}\,\int_0^1v_n''\,|v_n|^2\, dx=2\,|I|\;.
\end{equation}
Note that 
\[\int_0^1|v_n|^4\, dx=\int_0^1|u_n|^4\, dx=1\quad\forall\; n\in\N\;.
\]
Depending on the asymptotic behaviour of $(T_n)_{n\in\N}$, there are two possible cases~: 
\begin{description}
\item{(i)} If $\limsup_{n\to\infty}T_n=0$, then, because of (\ref{rel}), 
{$\limsup_{n\to\infty}\int_0^1v_n''{}^2\, dx\kern -0.5pt=\kern -0.5pt 0$}. Therefore, as $v'_n(0)=0$, $v'_n(x)=\int_0^x
v''_n(t)\,dt$, and $(v'_n)_{n\in\N}$ uniformly converges to $0$. Since $v_n$ cancels in $(0,1)$, the same
argument shows that $(v_n)_{n\in\N}$ uniformly converges to $0$. This is a
contradiction with the assumption that {$\int_0^1|v_n|^4\, dx=1$} for any $n\in\N$.
\item{(ii)} Up to the extraction of a subsequence, $(T_n)_{n\in\N}$ converges to some finite limit
$T$ in $(0,+\infty)$. Then $\int_0^1v_n''{}^2\, dx$ is uniformly bounded and, up to the
extraction of a further subsequence, $(v_n)_{n\in\N}$ weakly
converges in $H^2(0,1)$ and uniformly to some function~$v$ which is even, $1$-periodic and
non-increasing over the half-period for the same reason as in Proposition \ref{prop:reduc}. By Rellich's compactness theorem, $(v_n)_{n\in\N}$ strongly
converges in
$L^4(0,1)$ and {$\int_0^1 v^4dx\!=\!1$}. Due to~\eqref{rel2} and denoting $u(x):=T^{-1/4}\,v(x/T)$, we
get 
\begin{eqnarray}\label{liminf}
|I|&=&\liminf_{n\to\infty}T_n^{-7/2}\,\int_0^1v_n''{}^2\, dx\nonumber\\
&&\geq T^{-7/2}\,\int_0^1v''{}^2\,
dx=\int_0^T u''{}^2\,
dx
\end{eqnarray}
and
\begin{eqnarray}\label{lim}
2\,|I|&=&\lim_{n\to\infty}T_n^{-7/4}\,\int_0^1v_n''\,|v_n|^2\, dx\nonumber\\
&&=T^{-7/4}\,\int_0^1v''\,|v
|^2\, dx=\,\int_0^T u''\,|u
|^2\, dx\;,
\end{eqnarray}
together with $\int_0^T u^4\,dx=1$. Owing to the two facts that $u$ is $2\,T$ periodic and even, it is
then straightforward to check that $u$ is a minimizer for~$I_T$. As a consequence the inequality in~(\ref{liminf}) is an equality and up to the extraction of a subsequence $(v_n)_{n\in\N}$ strongly converges to $v$ in 
$H^2_{\rm loc}(\R)\cap C^{1,1/2}$. In particular $u'(0)=u'(T)=0$ holds and $u$ is non-increasing and changes sign on the half-period.
\end{description}
This ends the proof of the existence of a minimizer, after an eventual rescaling according to
Lemma~\ref{lem:scaling}. The Euler-Lagrange equation~(\ref{eq:EL-per}) is easily deduced, as already noted in Lemma \ref{Elr}. 

Moreover $u$ is decreasing on $(0,T)$. If it was not the case, $u$ would be constant on some interval and by the Euler-Lagrange equation this constant would be $0$. But since $u$ is not identically $0$ (because of $\int_0^T u^4\, dx=1$), this would be a contradiction with the Cauchy-Lispschitz theorem. 
\qed
\begin{remark}{\rm [Lower bound for the period]}\label{rem:period}
One can give an explicit lower bound for the value of $\liminf_{n\to\infty}T_n$ as follows~: 
\[
|v'_n(x)|^2=\left(\int_0^x v''_n(t)\,dt\right)^2\leq x\,\int_0^{1/2} |v''_n(t)|^2\,dt\;,
\]
if $0\leq x\leq 1/2$, whereas
\[
|v'_n(x)|^2=\left(\int_x^1 v''_n(t)\,dt\right)^2\leq (1-x)\,\int_{1/2}^1 |v''_n(t)|^2\,dt\;,
\]
if $ 1/2\leq x \leq1$. Thus 
$$\Vert v'_n\Vert_{L^\infty}\leq \sqrt{\frac{1}{2}\,\int_0^1 |v''_n|^2\,dx}$$ 
and since $v_n$ changes sign 
in $(0,1)$, 
\[
\Vert v_n\Vert_{L^\infty}\leq \Vert v'_n\Vert_{L^\infty}\;.
\] 
Thus, 
\[
1=\int_0^1 |v_n|^4\,dx\leq \frac 1 4\,\left(\int_0^1 |v''_n|^2\,dx\right)^2
\sim\frac 1 4\, T_n^7\,|I|^2\;,
\]
from which we deduce that
\[
\liminf_{n\to\infty}T_n \geq (|I|/2)^{-2/7}
\;.\]
\end{remark}
We are now going to prove that such a minimizer is unique. We begin with the following~:
\begin{lemma}{\rm [Infimum]}\label{Infimum}
Any periodic solution of 
\begin{equation}
\label{EL1}\left\{\begin{array}{l}
u^{(iv)}-2u\,u''-{u'}^2+2\,\lambda\,|u|^2\,u=0\\ \\
u(0)=1\;,\quad u''(0)=-a\;,\quad u'(0)=u'''(0)=0
\end{array}
\right.
\end{equation}
with $\lambda=-Q_{T(u)}(u)$ satisfies
\[-u''(0)=\sqrt{-Q_{T(u)}(u)}\; .\]
\end{lemma}
Recall that in the case of the Euler-Lagrange equations, $\lambda=-I$. \smallskip

\proof We denote $T=T(u)$ and $\lambda=-Q_{T(u)}(u)$ to lighten the notation. Multiply \eqn{EL1} by $u$
and $x\,u'$ and integrate on
$(0,2\,T)$~:
\begin{eqnarray*}
\int_0^{2\,T} {u''}^2\,dx+3\int_0^{2\,T} u\,{u'}^2\,dx+2\,\lambda\int_0^{2\,T} u^4\,dx =0\;,\\
T\,\Big(\lambda-|u''(0)|^2\Big)+\frac 32\int_0^{2\,T} {u''}^2\,dx+\int_0^{2\,T} u\,{u'}^2\,dx-\frac\lambda{2}\int_0^{2\,T} u^4\,dx =0\;.\end{eqnarray*}
Moreover, by definition of $\lambda$,
\[\int_0^{2\,T} {u''}^2\,dx+2\int_0^{2\,T} u\,{u'}^2\,dx+\lambda\int_0^{2\,T} u^4\,dx =0\;.\]
Combining these estimates, we get 
\[-Q_{T(u)}(u)=|u''(0)|^2\;,\]
which ends the proof.\qed
\begin{corollary}{\rm [Uniqueness]}\label{Uniqueness}
For a given period $T>0$, there is only one minimizer $u_T$ of $I_T$ which is even and
decreasing over the half period.
\end{corollary}
By the scaling invariance (Lemma \ref{lem:scaling}) we deduce that all such periodic minimizers $u_T$ are deduced from each other by a change of scale. 
\vskip6pt\noindent
\proof Uniqueness follows from Lemma~\ref{Infimum} if we prove first 
that $u'''(0)=0$. Assume that this is not the case and consider 
$\tilde u$ defined by~:
\[\tilde u(x)=\left\{\begin{array}{ll}
u(x)\quad &\mbox{if}\; x\in [0,T(u))\\ \\
u(-x)\quad &\mbox{if}\; x\in (-T(u),0)
\end{array}\right.\]
and extend it by periodicity. It is easy to check that
\[Q_{(-T(u),0)}(u)=Q_{(0,T(u))}(u)=I\;.\]
If this was not the case, say if $Q_{(-T(u),0)}(u)<Q_{(0,T(u))}(u)$, 
then we would indeed get
\[Q_{(-T(u),T(u))}(\tilde u)<I\;.\]
This means that $\tilde u$ is also a minimizer and solves the 
Euler-Lagrange equations on $(-T(u),$ $T(u))$. The function ${\tilde u}^{(iv)}$ is 
bounded in $L^2(-T(u),T(u))$, which implies that
$\tilde u'''$ is continuous at $x=0$~: then, by unique continuation, 
$-\tilde u'''(0_-)=\tilde u'''(0_+)=u'''(0)=0$.\qed

\begin{proposition}{\rm [$I$ is not a minimum]}\label{prop:infimum} The infimum $I$ is not achieved by a function in $L^4(\R)$.\end{proposition}
\proof Let us prove it by contradiction. Assume that $I$ has a minimizer $u$ in $L^4(\R)$. Because of the Euler-Lagrange equations, $u$ is smooth and has to decay to $0$ at infinity. Because of the uniqueness of the solutions of the Euler-Lagrange equations, $u$ cannot have compact support. The function~$u$ has infinitely many critical points, otherwise it is easy to define the tails of~$u$ as $u_{|(-\infty,\underline x)}$ and $u_{|(\overline x,+\infty)}$ where $\underline x$ and $\overline x$ are respectively the smallest and the largest critical points of $u$. The contribution of the tails to $Q_\R(u)=Q_{(-\infty,\underline x)}(u_{|(-\infty,\underline x)})+Q_{|(\underline x,\overline x)}(u)+Q_{(\overline x,+\infty)}(u_{|(\overline x,+\infty)})$ is clearly not optimal, for the same reason as in the proof of Proposition \ref{prop:ExistPeriodic} (case $T_n\to +\infty$).

Between two critical points, $u$ solve \eqn{Nehari} and is therefore made of the half of a periodic function. By Corollary \ref{Uniqueness}, the solution is uniquely determined, which means that $u$ itself is periodic. This is clearly a contradiction with the assumption that $u$ belongs to $L^4(\R)$.\qed

\section{A numerical computation of the infimum}\label{Numerics}

In this last section, we rely on the properties of the particular periodic minimizers
which have been built in the previous section -- that is, minimizers which are even with an absolute maximum at $0$ (up to a translation) and decrease over the half period -- to 
provide schemes in order to numerically compute the value of the infimum
$I$.\medskip

Any such minimizer of the periodic problem solves the
Euler-Lagrange equation~(\ref{eq:EL-per}) and satisfies 
\[u(0)=\max_{(-T(u),T(u))}\,u\;,\]
so that 
\[u'(0)=u'''(0)=0\;.\]
Furthermore, up to a rescaling (which means that one has to change the period accordingly), we may assume that 
\[u(0)=1\;.\]
This reduces the problem of finding a solution to a shooting problem in terms of 
\[a=-u''(0)\;,\]
once the value of $I$ is known. To determine $I$, we shall therefore proceed as follows~: Determine first the parameter $a\geq 0$
such that for $\lambda<1/4$, Equation~\eqn{EL1} has a solution $u$ such that $u'$ changes sign. It is an open question to determine theoretically the range of $\lambda$ and $a$ for which such a solution exists. Numerically, for $\lambda\in (9/64,1/4)$ we find
$a=a(\lambda)$ by a shooting method as follows~: For $a$ and $\lambda$ given, we solve \eqn{EL1} and
define $T(a,\lambda)$ as the first positive critical point of the solution, say $u_{a,\lambda}$~:
\[T(a,\lambda):=\inf\{x>0\;:\; u_{a,\lambda}'(x)=0\}\;.\]
It is not clear that such a quantity is always well defined and finite since $u_{a,\lambda}$ can be monotone decreasing or can even eventually explode monotonically. However, numerically we observe that this quantity makes sense. 

Then for a fixed $\lambda$ we minimize $Q_{(0,T(a,\lambda))}(u_{a,\lambda})$ on the set of the positive~$a$ for which $T(\lambda):=T(a,\lambda)$ is finite. This determines $a(\lambda)$. By periodicity, we extend the function $u_{a(\lambda),\lambda}$ from $(-T(\lambda),T(\lambda))$ to $\R$. Denote this extension by $u_\lambda$. There is no reason why $u_\lambda'''$ should be continuous at $x=T(\lambda)$, and in general $u_\lambda$ is not a solution of \eqn{EL1} on $\R$. Note that by construction $u'''_\lambda(0)=0$ and $u_\lambda$ is even. 
Then we minimize again 
\[J(\lambda):=Q_{(0,T(\lambda))}(u_{\lambda})\;.\]
It is easy to prove that if $J(\lambda)$ is well defined on $(9/64,1/4)$, then
\[\inf_{\lambda\in (9/64,1/4)}J(\lambda)=I=J(-I)\;.\]
Note that $u_I$ is a solution of \eqn{EL1} on $\R$. Numerically, we find \begin{center}\begin{tabular}{|c|}\hline\cr
$\displaystyle I\approx \,-\,0.1580\ldots$\cr\cr\hline \end{tabular}\end{center}
\begin{figure}[!ht]\begin{center}\includegraphics[scale=1]{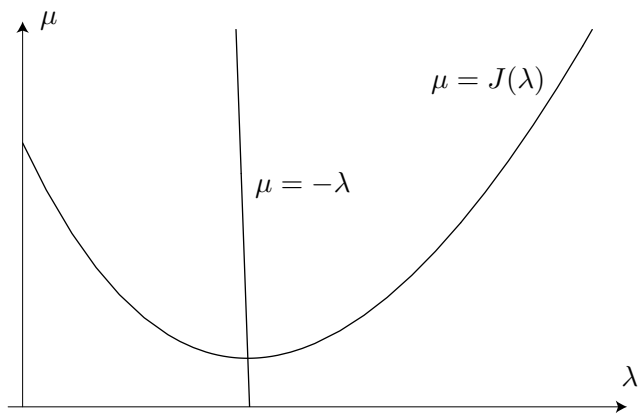} \caption{\small\sl Plot of the
function $J(\lambda)$ when $\lambda$ varies in the interval $(0.1,1/4)$. Note that $0.1<9/64\approx
0.140625\ldots<1/4=0.25$. The minimum $I=-\lambda$ is given as the solution of $J(\lambda)=-\lambda$. Also
note that the scales are not the same for $\lambda$ and $\mu$. }\label{Fig1}\end{center}\end{figure}
Alternatively, we can take advantage of the property stated in Lemma \ref{Infimum}. Define $\tilde u_\lambda$ 
as the solution of \eqn{EL1} with $a=\sqrt\lambda$ and
\[\tilde T(\lambda):=\inf\{x>0\;:\; \tilde u_\lambda'(x)=0\}\;.\]
There is again no {\sl a priori\/} reason why this quantity should be finite or even well defined, but numerically this makes sense. If we compute
\[\tilde J(\lambda):=Q_{(0,\tilde T(\lambda))}(\tilde u_{\lambda})\]
and numerically solve the equation
\be\label{eq:Alternative}\tilde J(\lambda)+\lambda=0\;,\ee
it is easy to check that
\[I=\tilde J(-I)\;.\]
In practice, the curves $\lambda\mapsto J(\lambda)$ and 
$\lambda\mapsto \tilde J(\lambda)$ are almost the same but the second method is much
more efficient. Numerically, Equation \eqn{eq:Alternative} has a single solution for $\lambda\in (0.1,1/4)$. Note that on $\R$, $\tilde u_{\lambda}$ is a solution of \eqn{EL1} which is not
necessarily $2\,\tilde T(\lambda)$ periodic. However, exactly as in the first method, it is a good family of test functions since it
contains the minimizer.
\begin{figure}
[!ht]
\begin{center}
\includegraphics[scale=1]{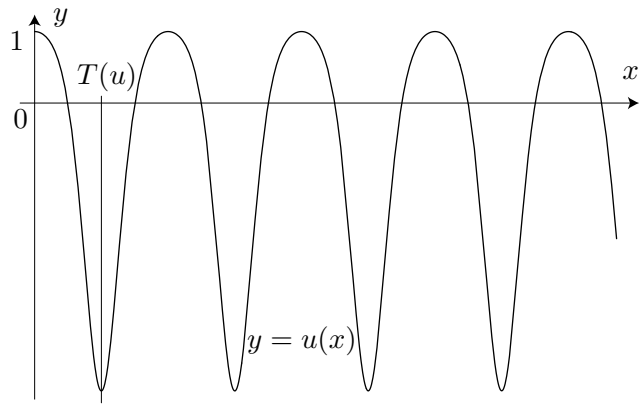} 
\caption{\small\sl Plot of the unique minimizer
$x\mapsto u(x)$ such that $u(0)=1$, $u'(0)=u'''(0)=0$, $u''(0)=-\sqrt{|I|}$. We note that $u$ changes sign and is monotone on $(0,T(u))$ with $T(u) \approx 3.43963\ldots$}
\label{Fig2}
\end{center}
\end{figure}

\medskip Let us conclude with a remark on quadratures. The equation for a minimizer can be reduced to a first order one as follows. Consider the Euler-Lagrange equation:
$$u^{(iv)}-2\,u\,u''-{u'}^2-2\,I\,u^3=0\;,$$
multiply it by $u'$ and integrate :
$$u'\,u^{(iii)}-\frac 12\,(u'')^2 -(u')^2\,u-\frac12\, I\,u^4=0\;,$$
where it is easy to check that the constant of integration is zero. Then we can set
$$u'=F(u)$$
for an unknown function $F$, and it is easy to check that the function {$y:=F^{3/2}$} satisfies:
$$y''=\frac32\, y^{-1/3}\, u + \frac34\, I\, u^4\, y^{-5/3}\;.$$
From the scaling property of the original equation, it follows that $y$ can be written as
$$y'=u^{5/4}\,f\left(\frac{y}{u^{9/4}}\right)\;,$$
where $f=f(z)$, $z=y/u^{9/4}$, has to satisfy
$$f'\,(f-\frac94\, z)=-\frac54\, f +\frac32\, z^{-1/3} +
\frac34\, I\, z^{-5/3}\;.$$
Thus our particular fourth order ODE can be reduced by successive quadratures to the integration of a first order ODE.



\end{document}